# Wind Turbine Design: Multi-Objective Optimization

Adam Chehouri, Rafic Younes, Adrian Ilinca and
Jean Perron

Additional information is available at the end of the chapter

http://dx.doi.org/10.5772/63481

**Abstract**

Within the last 20 years, wind turbines have reached matured and the growing worldwide wind energy market will allow further improvements. In the recent decades, the numbers of research papers that have applied optimization techniques in the attempt to obtain an optimal design have increased. The main target of manufacturers has been to minimize the cost of energy of wind turbines in order to compete with fossil-fuel sources. Therefore, it has been argued that it is more stimulating to evaluate the wind turbine design as an optimization problem consisting of more than one objective. Using multi-objective optimization algorithms, the designers are able to identify a trade-off curve called Pareto front that reveals the weaknesses, anomalies and rewards of certain targets. In this chapter, we present the fundamental principles of multi-objective optimization in wind turbine design and solve a classic multi-objective wind turbine optimization problem using a genetic algorithm.

## 1. Introduction

The exhaustion of fossil-fuel reserves, stricter environmental regulations and the world's rising energy needs have led to the deployment of renewable and sustainable energy sources. Among these alternatives, wind energy is a promising technology and recorded the fastest growing installed alternative-energy production according to reference [1]. It is expected that by the year 2030, at least 20% of the United States energy will be supplied by onshore and offshore wind farms [1]. In the next decade and half, it is vital that authorities record a significant increase in wind turbine installations and operability. Nevertheless, the prime conflict will continue to be the ongoing challenge to maintain a profitable and competitive cost of energy with fossil-fuel sources. Throughout the last 30 years, wind turbines have grown in size in order to reduce the cost of energy typically expressed in $/kWh. As a result, structural performance, durability requirements, safety hazards, transportation complications, noise and aesthetic pollutions all become issues that are more challenging for designers. Moreover, energy policies, international treaties, legislations and regulations set by governments have to be respected. For this reason, resolving the complex design problem of wind turbine design can be only achieved by *optimization* where an optimal solution is located. Many objective functions, design constraints, algorithms, tools and models have been proposed as will be discussed.

The rapid growth in the number of research papers on wind turbine design and optimization (**Figure 1**) during the last two decades highlights the status of the field of wind turbine optimization. In the past, some authors have compared the impact of different optimization objectives on the quality of the solution, others have reviewed the optimization algorithms, energy policies, economics, environmental impacts of wind turbines, but numerous researchers have proposed different optimization methodologies and resolution strategies. According to a study conducted by Chehouri et al. [2] in 2014, it was identified that less than 25% of the surveyed wind turbine optimization problems were solved using a multi-objective algorithm. In fact, solving such problems is not a straightforward task and often requires innovative techniques and algorithms. However, the main benefit is the ability to draw trade-off curves that reveal weaknesses and anomalies in the wind turbine design. For instance, it is very useful to pursue both minimum cost of energy and maximum annual energy production in the early stage of the design process.

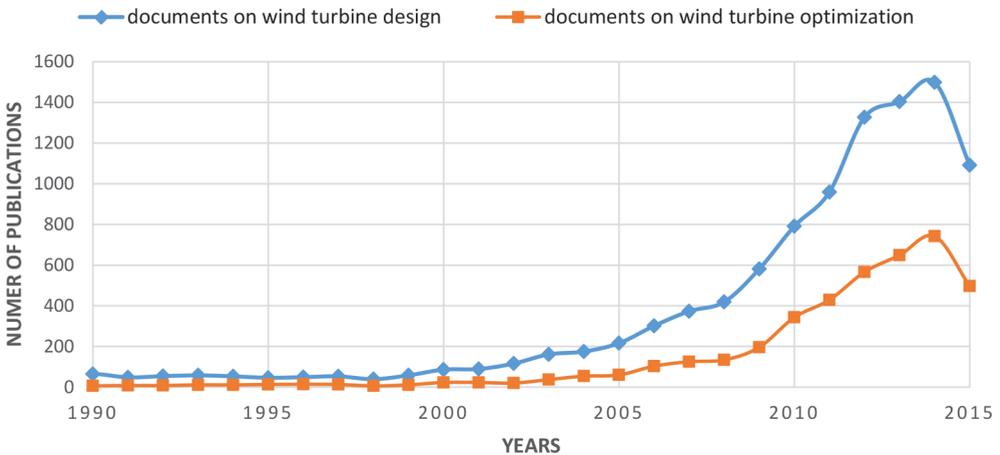

**Figure 1.** Number of published documents on wind turbine design in the last 40 years (reproduced from Scopus database).

Many authors carried a multidisciplinary study [3–12], where many objectives are considered in the design of wind turbines. The most common technique to combine conflicting functions (such as annual energy production and cost of energy) is by means of an appropriate set of weights. The variations that exist among these contradictory functions are essential for designers and therefore pursue to sketch the Pareto fronts.

In order to undertake the design of a horizontal wind turbine under multi-objective optimization (MOOP), there are numerous issues to be considered. The motif of this chapter is to present the fundamental principles of multi-objective optimization in the design of wind turbines. At the outset of this chapter in Section 2, we briefly discuss the fundamentals and terminology of multi-objective optimization. Section 3 highlights the objective functions that are used by designers. The design constraints that are enforced by wind turbine designers are enumerated in Section 4. We list the most relevant multi-objective optimization applied in wind turbine design in Section 5. The most common optimization algorithms used to solve multi-objective wind turbine optimization problems are examined in Section 6. Finally, a numerical example that demonstrates the resolution of a multi-objective design problem using a genetic algorithm (GA) is presented in Section 7.

## 2. Multi-objective optimization: state of art

### 2.1. Prologue

The term optimization refers to the finding of one or more feasible solutions, which correspond to extreme values of one or multiple objectives. Optimization methods are important in scientific experiments, particularly in engineering design and decision-making. When the problem is to find the optimal solution of one objective, the task is called *single-objective optimization*. There exist many algorithms that are gradient-based and heuristic-based that solves single-objective optimization problems. Beside deterministic search techniques, the field of optimization has evolved by the introduction of stochastic search algorithms that seek to find the global optimal solution with more ease. Among them, *evolutionary algorithms* (EAs) mimic nature's evolutionary principles and are now emerging as popular algorithms to solve complex optimization problems.

In engineering optimization, the designers are sometimes interested in finding one or more optimum solutions when dealing with two or more objective functions. This is known as *multi-objective optimization* and in fact, most real-world optimization problems involve multiple objectives. In this case, different solutions produce trade-offs or conflicting situations among the different objectives. Not enough emphasis is usually given to multi-objective optimization and there is a reasonable explanation for that. The majority of the methods that solve multi-objective optimization problems transform multiple objectives into a single function. Therefore, most of the efforts have been invested in improving the single-objective optimization algorithms. The studies concentrate on the conversion of multi-objective into single-objective problems, the convergence, constraint-handling approaches and the speed of these single-objective techniques.

Let us discuss the fundamental difference between single- and multi-objective optimization by taking two conflicting objective functions as examples. Obviously, each objective function possesses a unique and different optimal solution. For instance, if one is interested in buying a house, the decision-making has to take into consideration the cost and the comfort. If the buyer is willing to sacrifice comfort, they will get a house with the lowest price. However, if money is not an issue, then the buyer is able to afford a house with the best comfort. Between these two extremes, there exist many house choices at various costs and comfort. Now the big question is among these trade-offs, which solution is the best with respect to both objectives? Ironically, no house among the trade-off choices is the best with respect to both cost and comfort. Without any further information about these solutions (in our case example the houses), no solution from the set of trade-off can be said to be better than any other. This is the fundamental difference between a multi-objective and a single-objective optimization problem. From a practical standpoint, after a set of trade-off solutions are found, the user will use higher-level information to determine the convenient solution.

**2.2. Multi-objective evolutionary algorithms**

The classical way to solve multi-objective problems is to scalarize multiples objectives with a relative preference vector. Since only a single optimized solution can be found in one simulation, evolutionary algorithms shined as interesting methods to solve MOOP. The main reason is, unlike classical methods, EAs use a population of solutions in each iteration and therefore the outcome of an EA is a population of solutions. This ability to find multiple solutions in one single run made EAs an ideal approach to solve multi-objective optimization problems.

According to the available literature, the first real application of evolutionary algorithms in the determination of trade-off solutions for a MOOP was proposed in the doctoral dissertation of David Schaffer [13]. He developed the vector-evaluated genetic algorithm (VEGA), which demonstrated the ability of genetic algorithm to capture multiple trade-off solutions. Not much attention was given until another half a decade when David E. Goldberg published his book in 1989 [14] on a multi-objective evolutionary algorithm (MOEA) using the concept of dominance. From the latter derived many MOEAs such as Srinivas and Deb's non-dominated sorting (NSGA) [15] and the niched Pareto GA by Horn et al. [16]. Other techniques different from the domination-based MOEAs were proposed by Kursawe in 1990 (Kursawe's diploidy approach [17]) and Hajela and Lin's weighted-based approach [18] just to name a few. It can be easily seen from **Figure 2** that the number of studies conducted on the topic of MOEA has increased well over the last two decades. In less than 10 years, the number of year-wise papers has tripled and it can be expected that the growth will continue as new studies, books, surveys, research papers and dissertations will be published.

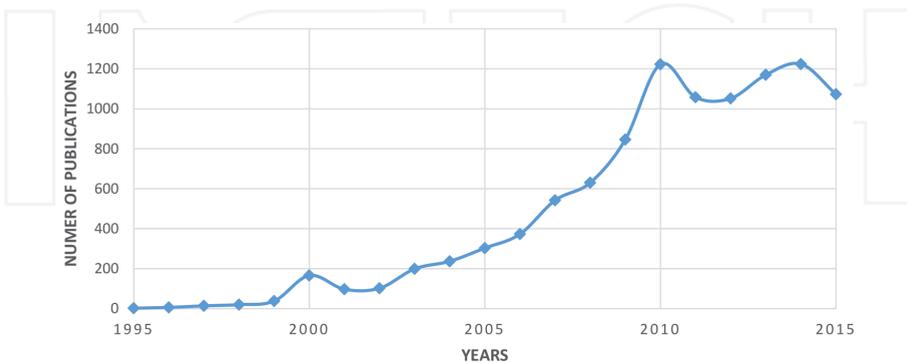

**Figure 2.** Number of published documents on multi-objective evolutionary algorithms.

**2.3. Multi-objective optimization problems**

A multi-objective optimization problem is composed of a number of objective functions, which are to be maximized or minimized. Similar to single-objective problems, the MOOP is subjected to a set of design constraints, which any optimal solution must satisfy. We can state the general form of a multi-objective optimization problem as follows:

$$\begin{aligned}&\textit{Minimize or Maximize } f_m(\vec{x}), m = 1, 2 \ldots M\\&\textit{Subject to } g_j(\vec{x}) \geq 0 \; j = 1, 2 \ldots J\\&h_k(\vec{x}) = 0 \; k = 1, 2 \ldots K\\&\vec{x}_i^L \leq \vec{x}_i \leq \vec{x}_i^U \; i = 1, 2 \ldots N\end{aligned} \qquad (1)$$

The solution $\vec{x}$ is a vector of $n$ variables $\vec{x} = (x_1, x_2, \ldots x_n)^T$. Often, the user will restrict the design variables between lower and upper bounds $\vec{x}_i^L$ and $\vec{x}_i^U$, respectively. In the above problem, there is $J$ inequality and $K$ equality constraints that can be linear and/or non-linear functions. A solution $\vec{x}$ is said to be *feasible* when *all* the constraints ($J + K + 2N$) are satisfied. Because of the presence of $M$ objective functions that need to be minimized and/or maximized, it is regularly convenient to apply the duality principle. The latter suggests that we can convert a maximization problem into a minimization one by multiplying the objective function by -1. This is a practical method because many optimization algorithms are developed to solve one type, for example, minimization problems. A major difficulty arises when any of the objective or constraint functions are non-linear, then the resulting MOOP becomes a non-linear multi-objective problem. Until now, the techniques to solve such problems do not have convergence proofs. Unfortunately, most real-world MOOPs are non-linear in nature and thus create a major challenge for scholars.

As stated earlier, the task in multi-objective optimization problems is to find a set of solution called the Pareto-optimal solution set, in which any two solutions must be non-dominated with respect to each other. In addition, any solution in the search space must be dominated by at least one point in the Pareto set. Therefore, the ultimate goal in multi-objection optimization is to find a set of solutions as close as possible to the Pareto-optimal front and as diverse as possible. The concept of domination is used in most MOOP algorithms. Without going into deep details, a solution $\vec{x}_1$ is said to dominate $\vec{x}_2$ if both conditions are satisfied:

1. The solution $\vec{x}_1$ is no worse than $\vec{x}_2$ in all objectives.

2. The solution $\vec{x}_1$ is strictly better than $\vec{x}_2$ in at least one objective.

To gain more knowledge on the procedures to find the non-dominated set in a given set *P* of size *N*, the reader is referred to [8].

## 3. Objective functions

The objective function that wind turbine engineers have used in their designs has evolved in the last few decades. In the early 1980s, the focus was towards the maximization of the power coefficient ($C_P$), which represents the theoretical fraction of power in the wind that can be extracted by the wind turbine. The blade shapes achieved with this strategy presented blades with large root chords and high twists. Moreover, since the maximization of the power coefficient occurs for a particular tip-speed ratio, the designer's interests shifted towards an alternative optimization metric—maximization of the annual energy production (AEP). Unlike the previous objective function, the AEP is obtained from a range of time (e.g. 1 year) and the given wind-speed spectrum of the wind turbine site. As wind energy is unable to compete with fossil fuel, the main objective has progressed to the minimization of the cost of energy (COE). The difficulty is the definition of the COE, mainly outlining the total costs of the wind turbine components. An earlier approach to decrease the COE is to limit the blade weight by restricting the chord lengths. Nevertheless, the stability of the blade is strongly affected by the decrease of the weight. Consequently, a proper balance between mass and stability must be ensured during the design of the wind turbine blade.

### 3.1. Maximization of the annual energy production

The annual energy is obtained by the integration of the wind turbine power curve with a wind-speed distribution (e.g. Weibull) over the wind-speed spectrum from the cut-in to cut-out speed [Eq. (2)]. One of the reasons why this metric is chosen by designers is due to the absence of a reliable structural and cost model. On the other hand, if the cost of energy is insignificant to the manufacturers and consumers, the maximum energy production is assumed as

$$AEP = \int_{V_{min}}^{V_{max}} P(V) f(V) dV \qquad (2)$$

$P(V)$ is the power curve of the wind turbine, $f(V)$ is the wind-speed distribution, $V_{min}$ is the cut-in speed and $V_{max}$ is the cut-out speed.

### 3.2. Minimization of the weight

Jureczko et al. [19] developed a numerical model of the wind turbine blade to perform a multi-criteria discrete-continuous optimization of wind turbine blades with the blade mass as the main objective function and the criteria's translated into constraints. Liao et al. [20] developed a multi-criteria-constrained design model with respect to minimum blade mass integrating a particle swarm optimization (PSO) algorithm using Federation Against Software Theft (FAST) [21]. Ning et al. [22] inspects the minimization of the turbine mass to AEP ratio as one of three examined objective functions. In a recent journal, Chen et al. [23] argue that a lighter blade mass will be beneficial to improve fatigue life based on requirements of the blade's strength and stiffness. Therefore, the minimum mass of the wind turbine blade was chosen as objective function.

## 4. Design constraints

Most engineering optimization problems include a set of equality and inequality constraints consisting of both linear and/or non-linear types. Generally, solving constrained optimization problems is more challenging than unconstrained systems. We have identified the most relevant constraints imposed by wind turbine blade designers as follows:

### 4.1. Ground clearance

A simple condition is set to prevent the collision of the blade with the ground.

### 4.2. Tip deflection

A constraint for the maximum tip deflection was included to ensure the local and global stability of the blade.

### 4.3. Shell and airfoil thickness

Some designers include a feasibility condition on the shell thickness and the surface of the airfoil to guarantee a proper trailing-edge separation.

### 4.4. Airfoil characteristics

In order to control the aerodynamic behaviour of the airfoil near stall, various constraints can be applied. For example, the absolute value of the slope beyond the stall angle can be regulated. Similarly, the coefficient of moment $C_{mc/4}$ and the ratio of the coefficient of lift and drag $C_L/C_D$ can be constrained to limit the blade torsion and undesirable separation behaviour, respectively. Another technique to avoid abrupt stall is to enforce a transition condition on each of the suction and pressure sides.

**4.5. Noise pollution**

A major objection for wind turbines is the noise that it generates as the blades rotate. Most of the aerodynamic noise models are semi-empirical and origin from the tip-vortex/trailing-edge interaction, turbulent inflow or the trailing-edge thickness.

**4.6. Stress**

The wind turbine is subjected to a large number of loads and therefore its components will be exposed to high stresses. To constrain these stresses, particularly on the wind blades and gearbox, the designers add inequality constraints that relate the generated stresses and the ultimate permissible stresses.

**4.7. Natural frequency**

In order to prevent the occurrence of resonance, the natural frequency of the blade must be separated from the rotor's rotation harmonics. Therefore, many designers limit some natural frequencies $\omega$ between an admissible bound [$\omega_{low}$; $\omega_{upper}$]. Another method is to assume a safety gap factor between the rotor's rated speed and the natural frequency. The reader is referred to references [24, 25] where the optimal design is pursued with respect to the maximum frequency design criterion.

# 5. Wind turbine MOOP

Researchers began examining multi-objective optimization algorithms in the design of wind turbines only two decades ago. In 1996, Selig and Coverstone-Carroll [5] examined the maximization of the AEP with no or few constraints on the loads. A short year later, Giguère and Selig [3] presented a multidisciplinary optimization (MDO) for the blade geometry of horizontal axis wind turbines (HAWTs). A sharing function [26] is used to obtain the trade-off curve between cost and energy. Only the structure of the blades is considered; however, the effects of the rotor on the wind turbine components are accounted in the cost model. The form of the cost model is indicated as follows (Eq. (3)):

$$C_i = C_{ib}\left(c_i + (1-c_i)\left(\frac{P_i}{P_{ib}}\right)\right) \quad (3)$$

where $C$ is the cost, $i$ is the component, $c$ is the cost factor (fixed portion of the total cost—cost factor of 20 % was used for the blades), $P$ is the design parameters of importance for component $i$ and subscripts ending with a $b$ are baseline values.

Benini et al. [6] apply an MOEA for the design optimization of stall-regulated wind HAWT with a trade-off between the ratio of AEP per wind park area (AEP$_{density}$ to maximize) and the cost of energy (minimize). The idea behind using the first metric is that the number of turbines that can be installed in a given area is inversely proportional to the square of the turbine radius. The MOEA handles the design parameters and searches for the optimal solutions following a set of Pareto concepts and basic principles of genetic programming [27, 28]. The authors choose the tip-speed, hub/tip ratio, and chord and twist distributions as their design variables. An assumption is made that the total turbine cost is reconstructed from the cost of the turbine blade alone, using the cost model from Eq. (4):

$$\text{CoE} = \frac{\text{TC} + \text{BOS}}{\text{AEP}} \text{FCR} + \text{OM} \tag{4}$$

where **TC** is the turbine cost, **BOS** the balance of station, **FCR** the fixed charge rate and **OM** the operation and maintenance costs.

In 2010, Grujicic et al. [10] developed a two-level optimization scheme to solve the MOOP. In the inner level, for a given aerodynamic design, the blade mass is minimized. In the outer level, a cost-assessment analysis is performed. This procedure is repeated until minimums are found for both the outer- and inner-level loops.

Kusiak et al. [9] introduced a data-driven approach to study the impact of turbine control on vibration. The authors developed a vibration prediction model using neural networks. Three objectives were included in the study (two vibrations and the power output), and a weighted sum of these objectives is minimized and set as follows (Eq. (5)):

$$\min\left(w_1 y_1(t) + w_2 y_2(t) + w_3 \frac{1}{y_3(t)}\right) \tag{5}$$

where $y_1(t)$ is the estimated vibration of the drive train, $y_2(t)$ the tower vibration model and $y_3(t)$ the estimated power output model.

In reference [29], Kusiak and Zheng present an approach to optimize the power factor and the power output of the wind turbine (bi-objective problem established by weights), using data-mining and evolutionary computation. The proposed approach generated optimized settings of the generator torque and the blade-pitch angle.

In the same year, Bottasso et al. [7] presented a thorough description of a multidisciplinary design optimization procedure. The authors assume that the weight is correlated to the cost. Since a reliable cost model is not offered to the public, no particular cost model was used. Instead of formulating a Pareto-optimal design problem, a combined cost was defined as the ratio of the annual energy production to the total weight. A two-stage sequential-constrained optimization algorithm was used to solve the constrained problem.

A year later, Wang et al. [11] presented a multi-objective algorithm. The power coefficient $C_P$ at the design wind speed of 9 m/s and the blade mass are chosen as the optimization objectives. The two objectives can be formulated as follows (Eq. (6)):

$$f_1 = \max(\boldsymbol{C}_p \,|\, V = 9\,\mathrm{m.s}^{-1})$$
$$f_2 = \min \int_{R_{hub}}^{R} \boldsymbol{m}_i dr \tag{6}$$

Where $C_p$ is the power coefficient, $m_i$ is the radial mass distribution of the wind turbine blade, and $R_{hub}$ and $R$ are the hub radius and the blade length, respectively.

In the design of HAWTs, the most obvious feature is the rapid growth in the size of HAWT blades [30]; therefore, it is insufficient to perform an airfoil shape optimization by itself. Rather, the entire blade geometry must be taken into consideration. According to reference [12], a wind turbine blade airfoil should satisfy the following aerodynamic requirements:

1. High lift-to-drag ratio ($C_L/C_D$).
2. High lift coefficient $C_L$.
3. Good performance during the stochastic behaviour of the wind flow.
4. Low sensitivity to leading-edge roughness.

The reader is referred to references [31–33] for the design of more efficient wind turbine blade airfoils. Ju et al. [12] developed a robust design optimization (RDO) for the design of a new series of wind turbine airfoils by maximizing both the $C_L/C_D$ and $C_L$. They also completed a sensitivity analysis of the roughness at the leading edge associated with the geometry profile uncertainty.

Multidisciplinary optimization was far and widely recognized of having the potential of becoming the cutting edge of the future [34, 35]. Kim [36, 37] and Michelena [38] apply maximum length sequence (MLS) algorithms called the target-cascading methodology method, where multiple levels and interfacing between the levels are defined.

For offshore wind turbines, the environmental conditions are more severe and more considerations have to be taken into account. Designers aim for a more proficient use of the capacity of the very expensive electrical cables, foundation, installation and erection costs.

## 6. Algorithms

### 6.1. Introduction

The selection of the appropriate optimization algorithm is a critical undertaking in any engineering optimization problem that relies on the attributes of the design space and on the nature of the problem. The final results depend on the algorithm in terms of accuracy and local minima sensitivity. Throughout the years, the algorithms used to solve wind turbine design problems have matured. At the outset, most of the methods derived directly from the blade-element momentum (BEM) theory, typically from Wilson and Lissaman [39]. In the 1990s, Selig and Coverstone-Carroll [5] were one of the originals to suggest a method based on GA for their wind turbine blade design tool. With the need to carry a multidisciplinary or multi-objective optimization design, Wood [40] and Sale et al. [41] simplified the MOOP into a single-objective question using a classical-weighted method. The approaches for solving conventional multi-objective design problems include:

- objective-weighted method
- hierarchical optimization method
- $\varepsilon$-constraint method
- goal-programming method.

It is important to highlight that all of the above algorithms convert the multi-objective problem into a single-objective problem. According to Ribeiro et al. [42], optimization algorithms can be categorized into two groups: gradient-based approaches (GBAs) and heuristic algorithms, whereas Endo [43] separates the optimization methods between genetic and non-genetic algorithms. In the last decades, in order to solve complicated optimization problems, evolutionary algorithms have been suggested such as:

- niched Pareto genetic algorithm (NPGA) [16]
- Pareto-archived evolution strategy (PAES) [44]
- strength Pareto evolutionary algorithm (SPEA) [45]
- neighbourhood cultivation genetic algorithm (NCGA) [46]
- non-dominated sorting genetic algorithm (NSGA)-II [47].

**6.2. Meta-heuristic algorithms**

Meta-heuristics are algorithms often inspired from nature, designed to replace or assist conventional solvers. This is a growing research field since the last few decades as meta-heuristics are now emerging as alternatives to the classical approaches.

An interesting fact can be drawn from the progress of the field of meta-heuristics and wind turbine optimization. In the recent years, wind energy showed an increase in the use of optimization methods such as linear programming, Lagrangian relaxation, quadratic programming and heuristic optimization (precisely genetic algorithm and particle swarm optimization) to name a few. However, it can be said that the gradient-based approaches and genetic algorithm are the two most popular optimization algorithms that have been applied in wind turbine design. The reason is quite simple, in the case of blade-geometry optimization, there is a large number of design variables, which are continuous (e.g. chord and twist distributions, blade pitch, etc.) and discrete (e.g. airfoil family, number of blades, etc.) at the same time. Moreover, some of these design variables are dependent from one another (e.g.

chord and twist), as well as competing objectives within the definition of the objective function (e.g. cost of energy).

GAs are the most popular evolutionary algorithms because of their robustness and reliability in wind turbine design problems. A genetic algorithm is an optimization method that mimics Darwin's principle of 'survival of the fittest' over a population of solutions (individuals) that evolves from one generation to another. It was originally proposed by Holland in 1975 [48]. Individuals with a large 'fitness' value have a superior probability to 'reproduce' in forming the new generation. Similar to a DNA chain, each individual is coded in one string and uses reproduction, crossover and mutation operators to direct the search over the generations. The usefulness of a GA is due to its robustness in multimodal design spaces. Likewise, GA explores non-linear, non-derivable, non-continuous domains and they are less sensitive to the initial domain.

## 7. Numerical example

In this section, we will solve a numerical example for the design of a wind turbine blade using a GA multi-objective optimization algorithm. The objective functions are the blade mass and the annual energy production. In order to calculate the mass, a structural model must be constructed. For the purpose of this study, a preliminary tool called Co-Blade [49] is used. As for the annual energy, WT-Perf [50] is introduced in the multi-objective platform to generate the AEP.

Co-Blade is a tool that helps designers to compute the structural properties of a wind turbine blade. It uses a combination of classical lamination theory (CLT) with an Euler-Bernoulli theory, and a shear-flow theory applied to composite beams is used to perform its analysis. This approach allows for a direct computation of the structural properties of a given blade, within several seconds of execution. The fitness function that Co-Blade minimizes is the blade mass penalized by the maximum stress, buckling, deflection and the natural frequency. The design variables are the chord-wise width of the spar cap at the inboard and outboard locations, the thickness of the 'blade-root' material and the thicknesses of the laminas within the leading-edge panel (LEP), trailing-edge panel (TEP), spar cap and shear webs along the length of the blade. They are listed in **Table 1**.

At first, the blade is represented as a cantilever beam under flap-wise and edge-wise bendings, axial deflection and elastic twist. Additional coupling between bending, extension and torsion is accounted for, due to the offsets between the beam-shear centre, tension centre and centre of mass from the blade-pitch axis (**Figure 3**). The beam cross sections are assumed to be thin-walled, closed and single- or multicellular, and the periphery of each beam cross section is discretized as a connection of flat composite laminates.

| Parameter | DESCRIPTION |
|---|---|
| w_cap_inb, w_cap_oub | Width of the spar cap normalized by the chord length at the INB_STN and OUB_STN blade stations |
| t_blade_root | Thickness of the 'blade-root' material at the INB_STN blade station |
| t_blade_skin1 …t_blade_skinN | Thickness of 'blade-shell' material at control points 1 through NUM_CP. The control points are equally spaced along the blade between the TRAN_STN and OUB_STN blade stations |
| t_cap_uni1 …t_cap_uniN | Thickness of 'spar-uni' material at control points 1 through NUM_CP |
| t_cap_core1 …t_cap_coreN | Thickness of 'spar-core' material at control points 1 through NUM_CP |
| t_lep_core1 …t_lep_coreN | Thickness of 'LEP-core' material at control points 1 through NUM_CP |
| t_tep_core1 …t_tep_coreN | Thickness of 'TEP-core' material at control points 1 through NUM_CP |
| t_web_skin1, t_web_skin2 | Thickness of 'web-shell' material at the two control points located at INB_STN and OUB_STN |
| t_web_core1, t_web_core2 | Thickness of 'web-core' material at the two control points located at INB_STN and OUB_STN. |
| w_cap_inb, w_cap_oub | Width of the spar cap normalized by the chord length at the INB_STN and OUB_STN blade stations |
| t_blade_root | Thickness of the 'blade-root' material at the INB_STN blade station |
| t_blade_skin1 …t_blade_skinN | Thickness of 'blade-shell' material at control points 1 through NUM_CP. The control points are equally spaced along the blade between the TRAN_STN and OUB_STN blade stations |
| t_cap_uni1 …t_cap_uniN | Thickness of 'spar-uni' material at control points 1 through NUM_CP |
| t_cap_core1 …t_cap_coreN | Thickness of 'spar-core' material at control points 1 through NUM_CP |
| t_lep_core1 …t_lep_coreN | Thickness of 'LEP-core' material at control points 1 through NUM_CP |

**Table 1.** Design variable for Co-Blade.

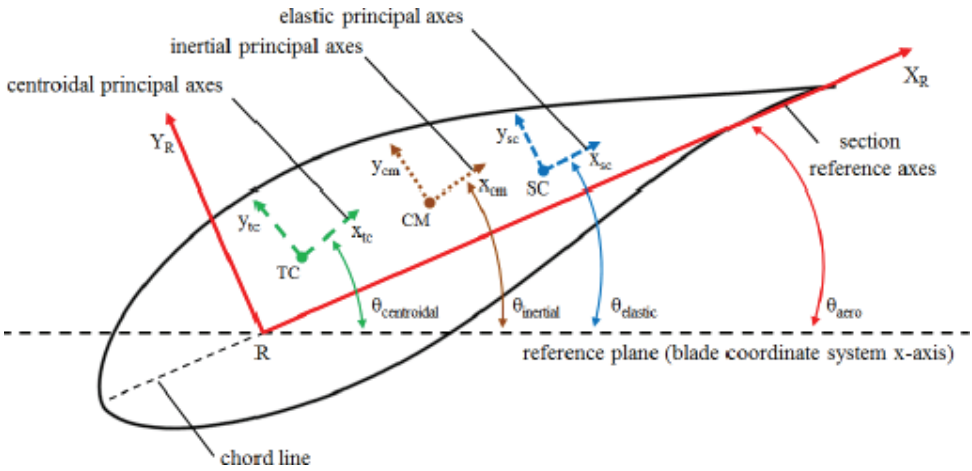

**Figure 3.** Orientation of the blade-axe systems [49].

In regard to Euler-Bernoulli beam theory, the beam cross sections are considered as heterogeneous and each of the material properties depends on the location in each cross section. The structural analysis at each discrete portion of the composite beam characterizes effective mechanical properties computed via classical lamination theory. Each discrete portion of the cross section then contributes to the global section properties of the composite beam (described further in references [51, 52]). Once the global cross-sectional properties are calculated, the deflections and effective beam axial stress ($\sigma_{zz}$) and the effective beam-shear stress ($\tau_{zs}$) can be now computed under the assumptions of an Euler-Bernoulli beam (refer to references [51–53]). The calculation of $\tau_{zs}$, the prediction of shear centre and torsional stiffness are based on a shear-flow approach, which is discussed in detail in reference [53]. Finally, by converting the distribution of effective beam stresses $\sigma_{ZZ}$ and $\tau_{zs}$ into equivalent in-plane loads, the lamina-level strains and stresses in the principal fibre directions ($\varepsilon_{11}$, $\varepsilon_{22}$, $\gamma_{12}$, $\sigma_{11}$, $\sigma_{22}$ and $\tau_{12}$) can be evaluated using CLT.

As mentioned earlier, Co-Blade applies a penalized blade mass defined as the following (Eq. (7)):

$$\text{Minimize } f(\bar{x}) = \text{Blade Mass } x \prod_{n=1}^{8} \max(1, p_n)^2 \qquad (7)$$

$$p_1 = \frac{\sigma_{11,max}}{\sigma_{11,fT}}$$

$$p_2 = \frac{\sigma_{11,min}}{\sigma_{11,fC}}$$

$$p_3 = \frac{\sigma_{22,max}}{\sigma_{22,fT}}$$

$$p_4 = \frac{\sigma_{22,min}}{\sigma_{22,fC}}$$

$$p_5 = \frac{|\tau_{12,max}|}{\tau_{12,y}}$$

$$p_6 = (\frac{\sigma}{\sigma_{buckle}})^\alpha + (\frac{\tau}{\tau_{buckle}})^\beta$$

$$p_7 = \frac{Tip\ Deflection}{Max\ Tip\ Deflection}$$

$$p_8 = max\left\{\frac{min\ allowable\ diff.between\ rotor\ freq.and\ the\ blade\ natural\ freq.}{|\omega_m - \omega_{rotor}|}\right\},$$

$$m = 1,\ldots Modes$$

Subject to

$\vec{x}_{LB} \leq \vec{x} \leq \vec{x}_{UB}$ (lower and upper bounds)

$A\vec{x} \leq \vec{b}$ (linear constraint)

Before we describe our fitness function, let us briefly discuss the second half of the multi-objective algorithm, the aerodynamic tool that calculates the AEP.

WT-Perf uses blade-element momentum theory to predict the performance of wind turbines with good accuracy. Users must build an appropriate input file that consists of the following set of data (**Table 2**):

- Model configuration
- WT-Perf algorithm configuration
- Cavitation model
- Turbine data
- Aerodynamic data
- Input/output settings
- Analysis settings.

| $E_{11}$ | $E_{22}$ | $G_{12}$ | $v_{12}$ | $\rho$ | Material Name |
|---|---|---|---|---|---|
| (Pa) | (Pa) | (Pa) | (-) | (kg/m³) | (-) |
| 2.80E + 10 | 1.40E + 10 | 7.00E + 09 | 0.4 | 1850 | (blade-root) |
| 2.80E + 10 | 1.40E + 10 | 7.00E + 09 | 0.4 | 1850 | (blade-shell) |
| 4.20E + 10 | 1.40E + 10 | 3.00E + 09 | 0.28 | 1920 | (spar-uni) |
| 2.60E + 08 | 2.60E + 08 | 2.00E + 07 | 0.3 | 200 | (spar-core) |
| 2.60E + 08 | 2.60E + 08 | 2.00E + 07 | 0.3 | 200 | (LEP-core) |
| 2.60E + 08 | 2.60E + 08 | 2.00E + 07 | 0.3 | 200 | (TEP-core) |
| 1.40E + 10 | 1.40E + 10 | 1.20E + 10 | 0.5 | 1780 | (web-shell) |
| 2.60E + 08 | 2.60E + 08 | 2.00E + 07 | 0.3 | 200 | (web-core) |

**Table 2.** Design variables for Co-Blade.

We have now defined two conflicting objective functions, the blade mass and the annual energy. Solving such MOOP can be achieved by the method of scalarizing. It consists of formulating a single-objective optimization problem such that optimal solutions to the single-objective optimization problem are Pareto-optimal solutions to the MOOP. A general formulation for a scalarization of a multi-objective optimization is given as (Eq. (8)):

$$\min \sum_{i=1}^{M} w_i f_i(\vec{x}) \qquad (8)$$

where the weights of the objectives $w_i > 0$ are the parameters of the scalarization.

We propose to use the following fitness function to minimize the mass and maximize the annual energy production (Eq. (9)):

$$\min\left(\alpha \frac{M(\vec{x})}{M_0} + (\alpha-1)\frac{\text{AEP}(\vec{x})}{\text{AEP}_0}\right) \qquad (9)$$

For a value of alpha near zero, the mass ratio is eliminated and the fitness function becomes $\min\left((\alpha-1)\frac{\text{AEP}(\vec{x})}{\text{AEP}_0}\right)$, which translates into the maximization of the normalized energy. Probably, for an alpha value close to 1, the energy ratio disappears and the problem is now a minimization of the mass. If we run the optimization problem for different values of alpha between 0 and 1, we can find Pareto-optimal solution to the MOOP. The reference mass and energy are taken, respectively, from the case study of alpha equals 0. The complete Pareto front is displayed in **Figure 4**.

Let us consider the following mechanical properties during the structural analysis. In our study, these properties are derived from Sandia 100-m blade SNL-100 [54]. **Table 3** lists the mechanical properties utilized in the structural design of the blade. Likewise, in **Table 2**, we list the configurations (input, model, turbine data and algorithm) for the input file required by the WT-Perf solver. The general flow chart of multi-objective optimization algorithm can be summarized in **Figure 5**. The complete inputs for the multi-objective optimization algorithm are listed in **Table 4**.

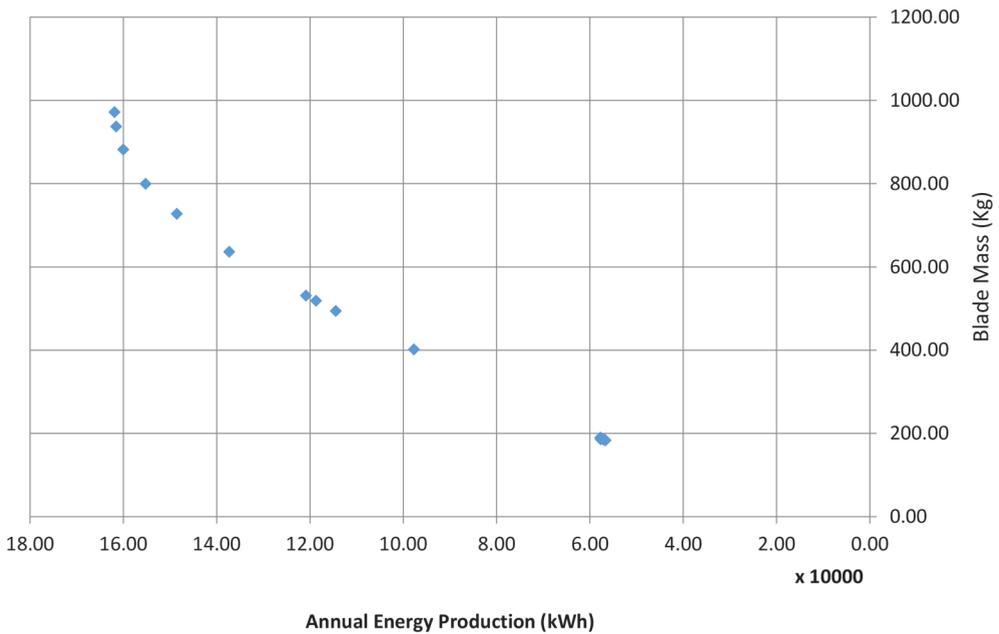

**Figure 4.** Pareto front for the given numerical example in Section 7.

| | | |
|---|---|---|
| **Input configuration** | | |
| False | Echo: | Echo input parameters to '<rootname>.ech'? |
| True | DimenInp: | Turbine parameters are dimensional? |
| True | Metric: | Turbine parameters are Metric (MKS vs FPS)? |
| **Model configuration** | | |
| 1 | NumSect: | Number of circumferential sectors. |
| 1000 | MaxIter: | Maximum number of iterations for induction factor. |
| 1.00E + 06 | ATol: | Error tolerance for induction iteration. |
| 1.00E + 06 | SWTol: | Error tolerance for skewed-wake iteration. |
| **Algorithm configuration** | | |
| True | TipLoss: | Use the Prandtl tip-loss model? |
| True | HubLoss: | Use the Prandtl hub-loss model? |
| True | Swirl: | Include Swirl effects? |
| True | SkewWake: | Apply skewed-wake correction? |
| True | AdvBrake: | Use the advanced brake-state model? |
| True | IndProp: | Use PROP-PC instead of PROPX induction algorithm? |
| True | AIDrag: | Use the drag term in the axial induction calculation? |
| True | TIDrag: | Use the drag term in the tangential induction calculation? |
| **Turbine data** | | |
| 3 | NumBlade: | Number of blades. |
| 10 | RotorRad: | Rotor radius (length). |
| 0.5 | HubRad: | Hub radius (length or div by radius). |
| 0 | PreCone: | Precone angle, positive downstream (deg). |
| 0 | Tilt: | Shaft tilt (deg). |
| 0 | Yaw: | Yaw error (deg). |
| 30 | HubHt: | Hub height (length or div by radius). |
| 30 | NumSeg: | Number of blade segments (entire rotor radius). |

**Table 3.** Input file for WT-Perf.

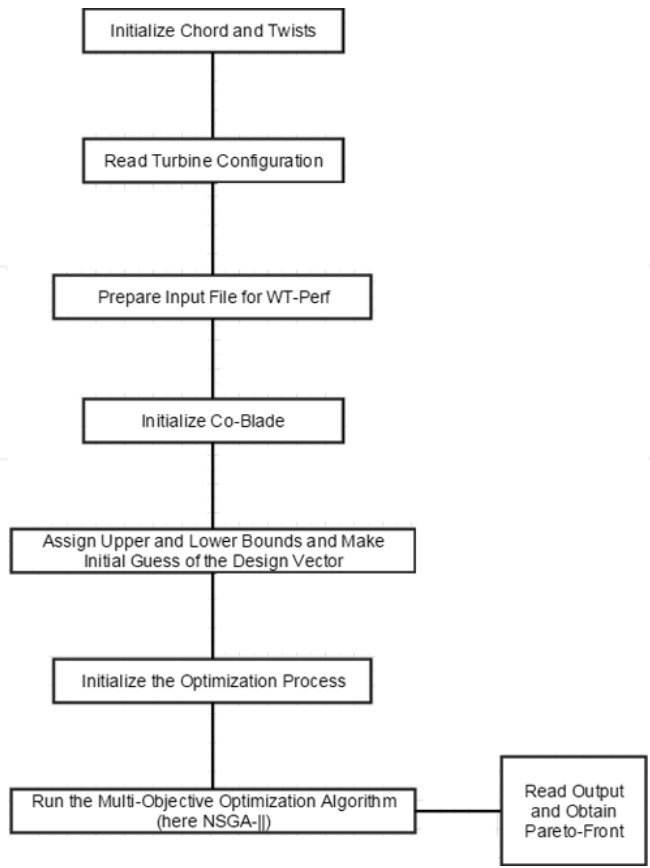

**Figure 5.** Flow chart of the multi-objective optimization algorithm.

| WT-Perf settings | | |
|---|---|---|
| 1000 | MaxIter: | Maximum number of iterations for induction factor. |
| 1.00E-06 | ATol: | Error tolerance for induction iteration. |
| 1.00E-06 | SWTol: | Error tolerance for skewed-wake iteration. |
| True | TipLoss: | Use the Prandtl tip-loss model? |
| True | HubLoss: | Use the Prandtl hub-loss model? |
| True | Swirl: | Include Swirl effects? |
| True | SkewWake: | Apply skewed-wake correction? |
| True | AdvBrake: | Use the advanced brake-state model? |
| True | IndProp: | Use PROP-PC instead of PROPX induction algorithm? |

| | **WT-Perf settings** | | |
|---|---|---|---|
| True | | AIDrag: | Use the drag term in the axial induction calculation? |
| True | | TIDrag: | Use the drag term in the tangential induction calculation? |
| 3 | | NumBlade: | Number of blades. |
| 0 | | Yaw: | Yaw error (deg). |
| 30 | | HubHt: | Hub height (length or div by radius). |
| 0.00001464 | | KinVisc: | Kinematic air viscosity |
| 0 | | ShearExp: | Wind-shear exponent (1/7 law = 0.143). |
| False | | UseCm: | Are Cm data included in the airfoil tables? |
| True | | TabDel: | Make output tab-delimited (fixed-width otherwise). |
| True | | KFact: | Output dimensional parameters in *K* (e.g. *kN* instead on *N*) |
| True | | WriteBED: | Write out blade-element data to '<rootname>.bed'? |
| True | | InputTSR: | Input speeds as TSRs? |
| 'mps' | | SpdUnits: | Wind-speed units (mps, fps, mph) |
| 0 | | NumCases: | Number of cases to run. Enter zero for parametric analysis. |
| WS or TSR | | RotSpd Pitch | Remove following block of lines if NumCases is zero. |
| 3 | | ParRow: | Row parameter (1-rpm, 2-pitch, 3-tsr/speed). |
| 1 | | ParCol: | Column parameter (1-rpm, 2-pitch, 3-tsr/speed). |
| 2 | | ParTab: | Table parameter (1-rpm, 2-pitch, 3-tsr/speed). |
| True | | OutPwr: | Request output of rotor power? |
| True | | OutCp: | Request output of Cp? |
| True | | OutTrq: | Request output of shaft torque? |
| True | | OutFlp: | Request output of flap-bending moment? |
| True | | OutThr: | Request output of rotor thrust? |
| 0.0 0.0 0.0 | | PitSt, PitEnd, PitDel: | First, last, delta blade pitch (deg). |
| 80 80 0.00 | | OmgSt, OmgEnd, OmgDel: | First, last, delta rotor speed (rpm). |
| **Analysis options** | | | |
| t | | SELF_WEIGHT: | Include self-weight as a body force? |
| t | | BUOYANCY: | Include buoyancy as a body force? |
| true | | CENTRIF: | Include centrifugal force as a body force? |
| true | | DISP_CF: | Apply correction factors to the beam displacements? |
| 0 | | N_MODES: | Number of modes to be computed |
| 50 | | N_ELEMS: | Number of blade finite elements to be used in the modal analysis |

| | **WT-Perf settings** | |
|---|---|---|
| **Optimization options** | | |
| t | OPTIMIZE: | Perform optimization of composite layup? |
| GS | OPT_METHOD: | Optimization algorithm for the optimization of composite layup |
| false | OPT_PITAXIS: | Optimize the pitch axis? |
| 0.375 | PITAXIS_VAL: | Pitch axis value outboard of max chord (ignored if OPT_PITAXIS = false) |
| 3 | INB_STN: | Inboard station where the leading- and trailing-edge panels, spar caps and shear webs begin |
| 8 | TRAN_STN: | Station where the root transition ends |
| 28 | OUB_STN: | Outboard station where the leading- and trailing-edge panels, spar caps and shear webs end |
| 4 | NUM_CP: | Number of control points between INB_STN and OUB_STN |
| false | READ_INITX: | Read the initial values for the design variables from INITX_FILE? |
| none | INITX_FILE: | Input file for the initial values of the design variables. |
| false | WRITE_STR: | Write structural input files at each function evaluation? |
| f | WRIT E_F_ALL: | Write the fitness value and penalty factors at each function evaluation? |
| f | WRIT E_X_ALL: | Write the design variables at each function evaluation? |
| f | WRITE_X_ITER: | Write the design variables at each iteration? |
| 100 | NumGens | Maximum number of generations for GA iterations |
| 100 | PopSize | Number of individuals per generation |
| 1 | EliteCount | Number of elite individuals per generation |
| 0.5 | CrossFrc | Fraction of individuals created by crossover |
| 1.00E-06 | GATol | Error tolerance for the GA fitness value |
| **Environmental data** | | |
| 1.225 | FLUID_DEN: | Fluid density (kg/m$^3$) |
| 9.81 | GRAV: | Gravitational acceleration (m/s$^2$) |
| 6.03 | U_mean: | Long-term mean flow (m/s) |
| 1.91 | Weib_k: | Shape factor |
| 6.8 | Weib_c: | Scale factor |
| **Blade data** | | |
| 30 | NUM_SEC: | Number of blade cross sections |
| 10 | BLD_LENGTH: | Blade length (m) |
| 0.5 | HUB_RAD: | Hub radius (m) |

| WT-Perf settings | | |
|---|---|---|
| 0 | SHAFT_TILT: | Shaft tilt angle (deg) |
| 0 | PRE_CONE: | Precone angle (deg) |
| 180 | AZIM: | Azimuth angle (deg) |
| 100 | MAX_ROT | Maximum rotational speed (rpm) |
| 10 | MIN_ROT | Minimum rotational speed (rpm) |
| cosine | INTERP_AF: | Interpolate airfoil coordinates? (choose "none", "cosine", or "equal" with no quotation marks) |
| 1 | ElmSpc | Blade-element radial spacing (0 equal, 1 cosinus) |
| 60 | N_AF: | Number of points in interpolated airfoil coordinates (ignored |
| mats-Wind.inp | MATS_FILE: | Input file for material properties |
| 0.13 | RootTranSt | Start of root transition region |
| 3 | RootTranSt_index | Index of start of root transition region |
| 0.288 | RootTranEnd | End of root transition region |
| 8 | RootTranEnd_index | Index of end of root transition region |
| 3 9 19 26 30 | CP_Index | Index of control points (chord and twist) |

**Table 4.** Input file for the multi-objective algorithm.

## 8. Conclusion

Within the last 20 years, wind energy conversion systems have reached maturity. The obvious growing worldwide wind energy market will culminate to further improvements. The continuous effort for the advancement in horizontal wind turbine performance strategies and techniques will result in additional cost reductions. The ultimate aim of any wind turbine manufacture is to design a wind turbine able to compete with fossil fuel. The number of research paper that applies optimization techniques in the attempt to reach an optimal blade design has demonstrated a significant increase in the recent decade alone. Despite the fact that a minimal cost of energy was chosen as the single main objective in most of the research papers, many have argued that it is more stimulating to evaluate the wind turbine design as an optimization problem consisting of more than one objective. Using multi-objective optimization algorithms, the designers are able to identify a trade-off curve called Pareto front that reveals the weaknesses, anomalies and rewards of certain targets. We can anticipate that future optimization problems will be set as multidisciplinary formulations. Consequently, solving such difficult optimization problem will require further developments in the optimization algorithm itself. Since traditional optimization techniques cannot overcome many of their drawbacks such as rapid divergence and sensitivity to the initial solution, population-based and nature-inspired algorithms will continue to emerge as worthy alternatives.

In this chapter, we presented the fundamental principles of multi-objective optimization in wind turbine design. We have identified the constraints and objective functions mostly targeted by designers. We briefly discussed the fundamentals and terminology of multi-objective optimization. The most common optimization algorithms used to solve multi-objective wind turbine optimization problems were presented. We highlighted the emergence of population-based techniques, particularly genetic algorithms. Finally, we showed the steps to solve a classic multi-objective wind turbine design problem using a genetic algorithm. The reader is referred to the following publications for further details [2, 55] concerning wind turbine optimization.

## Author details

Adam Chehouri[1,2*], Rafic Younes[2], Adrian Ilinca[3] and Jean Perron[1]

*Address all correspondence to: adam.chehouri1@uqac.ca

1 Anti Icing Materials International Laboratory (AMIL), Université du Québec à Chicoutimi, Québec, Canada

2 Faculty of Engineering, Third Branch, Lebanese University, Rafic Harriri Campus, Hadath, Beirut, Lebanon

3 Wind Energy Research Laboratory (WERL), Université du Québec à Rimouski, Québec, Canada